\documentstyle[twoside,fleqn,graphicx]{article}
\newcounter{zacountsec} 
\newcommand{\eh}{\hfill}\newlength{\sperr} 
\newcommand{\Title}[1]{{\large \bf #1}} 
\newcommand{\Author}[1]{{\sc #1}} 
\newcounter{symposium}
\newcommand{\Symposium}[1]{\setcounter{symposium}{#1}} 
\newcounter{session}
\newcommand{\Session}[1]{\setcounter{session}{#1}} 
\expandafter\ifx\csname jtitle \endcsname\relax\def\jtitle{}\fi 
\expandafter\ifx\csname symptitle \endcsname\relax\def\symptitle{}\fi 
\expandafter\ifx\csname makeheadings \endcsname\relax\def\makeheadings{}\fi 
\newcommand{\Section}[1]{{\stepcounter{zacountsec}\vspace{3mm}%
\hspace*{18mm}\normalsize\bf\arabic{zacountsec}. \parbox[t]{150mm}{ #1 }}} 
\newenvironment{Abstract}{\begin{minipage}[t]{177mm}\em }
{\end{minipage}} 
\newenvironment{thm}[2]{\begin{sloppypar}
{#1 #2.}\em{}}
{\end{sloppypar}}
\newcommand{\theo}[3]{\begin{thm}{#1}{#2} #3\end{thm}} 
\newcommand{\Theorem}{\hspace*{9mm}{\settowidth{\sperr}{\rm Theorem}%
\parbox[t]{1.3\sperr}{\rm T\eh h\eh e\eh o\eh r\eh e\eh m\eh } }}
\newcommand{\Lemma}{\hspace*{9mm}{\settowidth{\sperr}{\rm Lemma}%
\parbox[t]{1.3\sperr}{\rm L\eh e\eh m\eh m\eh a\eh } }}
\newcommand{\Proposition}{\hspace*{9mm}{\settowidth{\sperr}{\rm Proposition}
\parbox[t]{1.3\sperr}{\rm P\eh r\eh o\eh p\eh o\eh s\eh i\eh t\eh i\eh o\eh n\eh } }}
\newcommand{\Corollary}{\hspace*{9mm}{\settowidth{\sperr}{\rm Corollary}
\parbox[t]{1.3\sperr}{\rm C\eh o\eh r\eh o\eh l\eh l\eh a\eh r\eh y\eh } }}
\newcommand{\Definition}{\hspace*{9mm}{\settowidth{\sperr}{\rm Definition}
\parbox[t]{1.3\sperr}{\rm D\eh e\eh f\eh i\eh n\eh i\eh t\eh i\eh o\eh n\eh } }}
\newcommand{\proof}{\hspace*{9mm}{\settowidth{\sperr}{\rm Proof}%
\parbox[t]{1.3\sperr}{\rm P\eh r\eh o\eh o\eh f\eh. } }}
\newcommand{\N}{{\rm I}\!{\rm N}}
\newcommand{\R}{{\rm I}\!{\rm R}}
\newcounter{zalit}
\newenvironment{Acknowledgements}{\vspace{3mm}
\hspace*{18mm}{\bf Acknowledgements}\\[0.3cm]\begin{minipage}[t]{177mm}%
\small \em}{\end{minipage}}
\newenvironment{References}{
\Section{References}
\begin{small}\begin{list}{\arabic{zalit} }{\usecounter{zalit} 
\itemsep0mm \parsep0mm\settowidth{\labelwidth}{\small\rm 88}\labelsep0mm 
\setlength{\leftmargin}{\labelwidth}}}
{\end{list}\end{small}}
\newlength{\addrt}
\newenvironment{Address}{
\begin{minipage}[t]{\addrt}}
{\end{minipage}}
\newcommand*{\AuthorID}[1]{}\newcommand*{\LastName}[1]{}
\newcommand*{\FirstName}[1]{}\newcommand*{\ShortFirstName}[1]{}
\newcommand*{\Degree}[1]{}\newcommand*{\EMail}[1]{}
\newcommand*{\AuthorAddress}[1]{}
\newcommand{\AuthorInfo}[7]{%
\AuthorID{#1}\LastName{#2}\FirstName{#3}\ShortFirstName{#4}%
\Degree{#5}\EMail{#6}\AuthorAddress{#7}}

\begin{document}
\AuthorInfo{1}
    {Hirth}
    {Ulrich}
    {}
    {Dr.~rer.~nat.}
    {Ulrich.Hirth@UniBw-Muenchen.de}
    {Universit\"at der Bundeswehr M\"unchen, LRT 1, 
     Werner--Heisenberg--Weg 39, 85577 Neubiberg}
\AuthorInfo{2}
   {}
   {}
   {}
       {}
       {}
       {}
\Symposium{2}

\Session{0} 
\makeheadings \markboth{\jtitle}{\symptitle}
\newcommand{\veps}{\varepsilon}
\vspace*{15mm}\hspace*{18mm}
\begin{minipage}[t]{157mm}
\Author{Hirth, Ulrich}
\vspace*{0.4cm}

\Title{Extreme exchangeable random order processes
by positive definite functions on semigroups}

\end{minipage}
\vspace*{3.5mm}

\begin{Abstract}
\end{Abstract}

\Section{Motivation}

Often in practice jobs are partially ordered in such a way that the second job
cannot begin before the first one is finished.  Such situations are considered
by Kyle Siegrist in his recent paper [S]~.  So we have a partially ordered 
network of jobs and are interested in, among others, the subset of jobs which
are already completed at time $t$~.  Clearly this increases with $t$~, and the
mentioned subset can be viewed as partially ordered again, the order being
inherited by the original order.  So we have an increasing family of partially
ordered sets.  Such a situation asks itself for more algebraic structurisation.
In fact there is the following new application of the [HR]--method.

\Section{The set--up of the Bauer simplex of continuous exchangeable
  probability measures on the set of order processes}

We use here the notation and terminology of [HR2000]~.
In particular ${\cal V}$ is the set of all reflexive, transitive {\it (but not
necessarily anti-symmetric)} relations 
--- called by us {\bf partial orders} --- on $M$~, where $M$
is a countable infinite set, and we recall that ${\cal V}$ carries a natural
metrisable topology and that the diagonal relation $D$ is the neutral element
of the semigroup $({\cal V},\vee)$~.
\theo{\Definition}{1}{We denote by ${\cal Y}$ the set of all left--continuous
increasing maps $Y:[0,\infty[\ \mapsto{\cal V}$ with $Y(0)=D$~.}
\theo{\Definition}{2}{For $Y,Z\in{\cal Y}$ the join $Y\vee Z\in{\cal Y}$ 
is given by $(Y\vee Z)(t):=Y(t)\vee Z(t)$~.}
\theo{\Proposition}{3}{$({\cal Y},\vee)$ is an idempotent abelian semigroup
  with neutral element the constant $D$ and absorbing element the map
  which switches from $D$ to $M\times M$ immediately at $t=0$~.}
\theo{\Definition}{4}{By the {\bf support} of a $Y\in{\cal Y}$ we understand
  the set 
$\langle Y\rangle:=\bigcup_{t\in[0,\infty[} \langle Y(t)\rangle\subseteq M$~.}
\theo{\Definition}{5}{
${\cal Z}:=\{Y\in{\cal Y}\vert\langle Y\rangle\ \mbox{is finite.}\}$}
\theo{\Lemma}{6}{$\langle D\rangle=\emptyset$ and 
$\langle Y\vee Z\rangle\subseteq\langle Y\rangle\cup\langle Z\rangle$ for all
$Y,Z\in{\cal Y}$~, so ${\cal Z}$ is a sub--semigroup of ${\cal Y}$~.}
\theo{\Definition}{7}{For $Y,Z\in{\cal Y}$ we say $Y\le Z$ if and only if
$Y\vee Z=Z$~.}
\theo{\Lemma}{8}{This is the case if and only if $Y(t)\le Z(t)$ for all 
$t\in[0,\infty[$~.}
\theo{\Lemma}{9}{``$\le$'' is a partial order on ${\cal Y}$~.}
\theo{\Definition}{10}{For $Z\in{\cal Z}$ we put 
$Q_Z:=\{Y\in{\cal Y}\vert Z\le Y\}$~.}
\theo{\Lemma}{11}{For all $Y,Z\in{\cal Z}$ we have 
$Q_{Y\vee Z}=Q_Y\cap Q_Z$~.}
\theo{\Definition}{12}{A sub--semigroup $I$ of ${\cal Z}$ is called
  {\bf left--hereditary} if for all $Z_1,Z_2\in{\cal Z} \quad Z_1\le Z_2\in I$
  implies $Z_1\in I$~.}
\theo{\Proposition}{13}
{The join $\bigvee_{Z\in I} Z$ is left--continuous again.}
\begin{proof}
For this let us remember why the semigroup operation 
$\vee:{\cal V}\times{\cal V}$ is continuous.  We see ${\cal V}$ as a closed
subset of the compact metrisable space $\{0,1\}^{M\times M}$~.
There is a natural mapping $p:\{0,1\}^{M\times M}\mapsto{\cal V}$ which
assigns to any $\tilde V\in\{0,1\}^{M\times M}$ the smallest partial order
containing $\tilde V$~.  It is easily seen that $p$ is continuous, on the
basis of the explicit construction with the chains and because the subsets
$Q_U$ (~with $U\in{\cal U}$~) generate the topology of ${\cal V}$~.
The semigroup operation $\vee$ on ${\cal V}$ is now the obvious composition of
continuous maps $${\cal V}\times{\cal V}
\mapsto\{0,1\}^{M\times M}\times\{0,1\}^{M\times M}
\mapsto\{0,1\}^{M\times M}\mapsto{\cal V}\ ,$$
so $\vee:{\cal V}\times{\cal V}\mapsto{\cal V}$ is continuous.
The problem with infinitely many factors instead of just two factors is
the map in the middle in the above chain composition.  I do not see that it is
continuous.  However, for our left--continuity we only need continuity 
from below.  So, in terms of subsets of $M\times M$~, the question is, 
if $A_{mn}\nearrow A_m$ as $n\to\infty$ for each $m$~, whether or not
$\bigcup_m A_{mn}$ converges to $\bigcup_m A_m$~.  It is definitively
increasing in $n$~, so the question is whether or not
$$\bigcup_n \bigcup_m A_{mn} = \bigcup_m \bigcup_n A_{mn}$$ is true.
Trivially, it {\bf is} true, completing the proof of the proposition.
\end{proof}
\theo{\Proposition}{14}{If $I$ is a left--hereditary sub--semigroup of 
${\cal Z}$ and $Z\in{\cal Z}$ is such that $Z\le\bigvee_{Y\in I} Y$~,
then for each $\veps>0$ \quad $Z(t)$ is $\le$ the maximum of finitely many 
$Y\in I$~, taken at the time $t+\veps$~, for all $t$~.}
\begin{proof}
For all $j,k\in M$ we put
$t_0=t_0(j,k):=\inf\{t\in[0,\infty[\vert(j,k)\in Z(t)\}$ 
with the usual convention
$\inf(\emptyset)=\infty$~.
For each pair $(j,k)$ with $t_0(j,k)<\infty$ and each $\veps>0$
there is a $\tilde Y_{jk}\in I$
such that we have $(j,k)\in\tilde Y_{jk}(t_0+\veps)$~.
We define $Y_{jk}(t)$ to be $D$ on $[0,t_0+\veps]$ 
and to be $\tilde Y_{jk}(t_0+\veps)$ on $]t_0+\veps,\infty[$~, 
and because of left-hereditarity we have still 
$Y_{jk}\in I$~.  For $j,k\in\langle Z\rangle$ with $t_0(j,k)=\infty$ we put
$Y_{jk}$ to be the constant $D$~, which is always $\in I$~.
So we have $Z(t)\le\bigvee_{j,k\in\langle Z\rangle} Y_{jk}(t+\veps)$
with $\bigvee_{j,k\in\langle Z\rangle} Y_{jk}\in I$~, 
and the proposition is proved.
\end{proof}

\vspace{0.5cm}

Next, in order to endow ${\cal Y}$ with a suitable topology, so to introduce
in a useful way the space $M^1_+({\cal Y})$ of (Radon) probability measures on
${\cal Y}$~, we need to resort to a {\bf Skorokhod}--like topology on 
${\cal Y}$~.  For this we recall that ${\cal V}$ can be seen as a subset of
$\{0,1\}^{M\times M}$~, and we start with the basic set $\{0,1\}$ in the place
of ${\cal V}$~, denoting by ${\cal Y}^\prime$ the set of all left--continuous
increasing maps $Y^\prime:[0,\infty[\ \mapsto\{0,1\}$ 
with $Y^\prime(0)=0$~.
Any such map can be identified with its jump point 
$\in [0,\infty]$~, and we let ${\cal Y}^\prime$ inherit its
(compact) topology from there.
Our multidimensional set ${\cal Y}$ consists of certain maps from $[0,\infty[$
to ${\cal V}\subset\{0,1\}^{M\times M}$~, and viewed as maps from $[0,\infty[$
to $\{0,1\}^{M\times M}$ this splits into the component maps from $[0,\infty[$
to $\{0,1\}$~, which are all left--continuous and increasing.  In this way 

\vspace{0.5cm}

\theo{\Definition}{15}{
$${\cal Y}\mbox{ inherits its topology from the isomorphic set }$$
$$\tilde{\cal Y}:=\{f\in [0,\infty]^{(M\times M)\setminus D}
\ \vert\ f(j,l)\le\max\left(f(j,k)\ ,\ f(k,l)\right)
\ \mbox{for all pairwise different}\ j,k,l\in M\}\ ,$$

\vspace{0.5cm}

$[0,\infty]^{(M\times M)\setminus D}$ being endowed with the
product topology, which is compact by Tychonov's theorem and metrisable,
and ${\cal Y}$ being compact as a closed subset of 
$[0,\infty]^{(M\times M)\setminus D}$~.}
\theo{\Definition}{16}{By $M^1_+({\cal Y})$ we denote the space of all Borel
probability measures on ${\cal Y}$~, all of which are automatically Radon,
${\cal Y}$ being endowed with the above compact metrisable topology.}
\theo{\Proposition}{17}{The Borel $\sigma$--algebra on ${\cal Y}$ is generated
by the closed subsets $Q_Z$~.}
\begin{proof}
In the framework of the isomorphic description of ${\cal Y}$ given by
Definition~15~, denoting the canonical projections from $\tilde{\cal Y}$
to $[0,\infty]$ by $pr_{jk}$~, \quad $Q_Z$ corresponds to the intersection
of finitely many sets of the form $pr_{jk}^{-1}\left([0,t]\right)$~,
where the $t$'s are the switching times of $Z$ (~since we require the switching
times of the $Y$'s in $Q_Z$ to be $\le$ those of $Z$~), so the sets
$Q_Z$ are closed.\\
The Borel $\sigma$--algebra on $[0,\infty]^{(M\times M)\setminus D}$
is generated by the closed subsets, 
and the ``finite--dimensional'' closed subsets
(~``$\times [0,\infty]\times [0,\infty]\times\dots$''~) already suffice
to generate this, since any other closed subset is the intersection
of countably many ``finite--dimensional'' ones, 
the finite--dimensional projections being continuous.
As we know from $[0,\infty]^n$~,
the $\sigma$--algebra generated by the finite--dimensional closed intervals
$[\vec{\bf 0},\vec{\bf t}]$ is equal to the one generated by all closed sets,
and the intervals $[\vec{\bf 0},\vec{\bf t}]$ correspond to our sets $Q_Z$~.
Finally, the transition from $[0,\infty]^{(M\times M)\setminus D}$ 
to ${\cal Y}$ is achieved by taking the trace $\sigma$--algebra on the
closed subset $\tilde{\cal Y}$~, concluding the proof of the Proposition.
\end{proof}
\theo{\Theorem}{18}{A function $\varphi:{\cal Z}\mapsto\R$ is positive
definite with respect to $\vee$~, normalised (i.e.~$\varphi(D)=1$~) 
and continuous from below if and
only if $$\varphi(Z)=\mu(Q_Z)\quad\mbox{for all }Z\in{\cal Z}$$
for a (uniquely determined) $\mu\in M^1_+({\cal Y})$~.}
\begin{proof}
For any $\mu\in M^1_+({\cal Y})$ the function
$\varphi:{\cal Z}\mapsto [0,1]$ defined by $\varphi(Z)=\mu(Q_Z)$
is positive definite and fulfills $\varphi(D)=1$~, as is easily seen 
in the same way as in previous work of Paul Ressel and myself.
Moreover it is continuous from below, because this is totally
analogous to the continuity {\it from above} for distribution functions
of probability measures on (certain Borel--measurable subsets of)
$[0,\infty]^n$~, \quad $Q_Z=\{Y\in{\cal Y}\vert Y\ge Z\}$ being the event
that the switching times of $Y$ are all $\le$ those of $Z$~,
where only finitely many switching times of $Z$ are allowed to be $<\infty$~.\\
Conversely, if $\varphi$ is any positive definite function 
with $\varphi(D)=1$~, there is a unique Radon probability measure $\nu$ on
$\hat{\cal Z}$ representing $\varphi$ via
$$\varphi(Z)=\int 1_I(Z)\;d\nu(I)=\nu(\{I\in{\cal I}\vert Z\in I\})\ ,$$
where ${\cal I}$ denotes the set of all left--hereditary sub--semigroups of
${\cal Z}$~, and where the identification of $I\in{\cal I}$ 
with $1_I\in\{0,1\}^{\cal Z}$ is used to topologise ${\cal I}$~.
We define $h:{\cal I}\mapsto{\cal Y}$ by
$$h(I):=\bigvee_{Z\in I} Z \quad .$$
In the following argument we mean by the sub--index ``$-\veps$''
the shift operation $f\mapsto f(\bullet-\veps)$
(~where for times $t\le 0$ \quad $f(t):\equiv D$~)
and by the sub--index ``$+\veps$'' we mean
$f\mapsto f(\bullet+\veps)$~.
Note that in ${\cal Z}$ \quad $Z^{(1)}=Z^{(2)}$
is equivalent to $Z^{(1)}_{-\veps}=Z^{(2)}_{-\veps}$~,
but not equivalent to $Z^{(1)}_{+\veps}=Z^{(2)}_{+\veps}$~.
So in the argument below
$Z_{-\veps}\in I$ clearly {\it implies}
$Z=\left(Z_{-\veps}\right)_{+\veps}\in I_{+\veps}$~,
and $Z\in I_{+\veps}$ implies
$Z_{-\veps}\in \left(I_{+\veps}\right)_{-\veps}$~.
Here, for $I$~, first taking ``$+\veps$'' and afterwards
``$-\veps$'' means replacing the elements of $I$
on $[0,\veps]$ by $D$~, and the result of this is still
in $I$ because of the left--hereditarity of $I$~,
so we luckily arrive at $Z_{-\veps}\in I$ as desired.
With the first of the following ``$\supseteq$''
being by Proposition 14, we have
$$\{I\in{\cal I}\vert Z_{-\veps}\in I\}
=\{I\in{\cal I}\vert Z\in I_{+\veps}\}
\supseteq\{I\in{\cal I}\vert Z\le h(I)\}
=h^{-1}(Q_Z)
\supseteq\{I\in{\cal I}\vert Z\in I\}\ ,$$
where the $\nu$--measure of the leftmost resp.~rightmost side 
is $\varphi(Z_{-\veps})$ resp.~$\varphi(Z)$~.
For $\veps\searrow 0$ the switching times of $Z_{-\veps}$ 
converge from above to those of $Z$~, 
which means $Z_{-\veps}\nearrow Z$ in ${\cal Z}$~.
So with the additional assumption that $\varphi$ is continuous from below,
$\varphi(Z_{-\veps})$ converges to $\varphi(Z)$~.
This shows that $h$ is measurable with respect 
to the $\nu$--completion of the $\sigma$--algebra on ${\cal I}$
and that $\varphi(Z)=\nu(h^{-1}(Q_Z))=:\mu(Q_Z)$
with $\mu$ being the image measure $\nu^h\in M^1_+({\cal Y})$ as desired,
finishing the proof of the theorem.
\end{proof}
\theo{\Corollary}{19}{In $M^1_+({\cal Y})$ a sequence $(\mu_n)_{n\in\N}$
converges to $\mu$ if and only if the corresponding positive definite
functions fulfill $\varphi_n(Z)\to\varphi(Z)$ for all $Z\in{\cal Z}$ with
$\mu(\partial Q_Z)=0$~.}
\begin{proof}
By the portmanteau theorem, the implication ``$\Rightarrow$'' is obvious.
Conversely, if $\varphi_n(Z)\to\varphi(Z)$ for all $Z\in{\cal Z}$ with
$\mu(\partial Q_Z)=0$~, then by the compactness of $M^1_+({\cal Y}$ we have
along some sub--sequence $\mu_n\to\tilde\mu$ 
and $\varphi_n(Z)\to\tilde\varphi(Z)$ 
for all $Z$ with $\tilde\mu(\partial Q_Z)=0$~,
so $\varphi(Z)=\tilde\varphi(Z)$ 
for all $Z$ with $\mu(\partial Q_Z)=\tilde\mu(\partial Q_Z)=0$~.
By an analoguous argument as with distribution functions on $\R^n$~,
this suffices to imply $\mu=\tilde\mu$ and finally $\mu_n\to\mu$
for the whole sequence.
\end{proof}

Next we introduce our notion of {\bf exchangeability} for measures in
$M^1_+({\cal Y})$~.  The {\it permutations} (finite or infinite ones, this
does not matter here) of the countably infinite base set $M$ act in the
natural way on the set ${\cal V}$ of partial orders, cf.~[HR2000]~.
And on ${\cal Y}$ such a permutation acts timepointwise in the same way as it
acts on ${\cal V}$~, and so the action carries in the canonical way over
to $M^1_+({\cal Y})$~.  
\theo{\Definition}{20}{If a $\mu\in M^1_+({\cal Y})$ is invariant
under all these permutations, then we call it {\bf exchangeable} and denote
this in signs by $\mu\in M^{1,e}_+({\cal Y})$~.}

Next we express this property by means of a mapping $$g:{\cal Z}\mapsto T\ ,$$
where $(T,+)$ is another abelian (but not idempotent) semigroup.
$T$ is defined as the set of all isomorphy classes of ${\cal Z}$~,
where $Z_1,Z_2\in{\cal Z}$ are called isomorphic if they are
permutations of each other.  $g:{\cal Z}\mapsto T$ is then the mapping which
assigns to every $Z\in{\cal Z}$ its isomorphy class $g(Z)\in T$~.
The semigroup operation ``$+$'' on $T$ is defined by
$$g(Z_1)+g(Z_2):=g(Z_1\vee Z_2)\mbox{ for representants }Z_1\mbox{ and }Z_2
\mbox{ with disjoint supports.}$$
Clearly this is well--defined.  In the same way as in our previous papers it
is straightforward that $g$ is strongly almost additive and hence
{\bf strongly positivity forcing.}  The mapping $g$ now expresses
exchangeability as follows: let $\mu\in M^1_+({\cal Y})$~; then we have
the corresponding positive definite function $\varphi:{\cal Z}\mapsto\R$ with
$\varphi(Z)=\mu(Q_Z)$~, and $\mu$ is exchangeable if and only if $\varphi(Z)$
depends only on the isomorphy class of $Z$~, factorising over $g$~, 
i.e.~$\varphi=f\circ g$ for some mapping $f:T\mapsto\R$~.
This mapping $f$ is necessarily positive definite again, since $g$ is
positivity forcing.\\
Because $g$ is even {\it strongly} positivity forcing, the set 
${\cal P}^{1,g}({\cal Z})$ of positive
definite functions on ${\cal Z}$ factorising over $g$ is a {\bf Bauer simplex}
whose extreme points are precisely the functions of the form $\rho\circ g$
with $\rho$ being a bounded {\it character} on the semigroup $(T,+)$~.

\theo{\Definition}{21}{Let ${\cal P}^{c,g}({\cal Z})$ denote the set of those
functions in ${\cal P}^{1,g}({\cal Z})$ which are continuous from below.}
\theo{\Lemma}{22}{${\cal P}^{c,g}({\cal Z})$ is an {\bf extreme} subset
of ${\cal P}^{1,g}({\cal Z})$~, 
i.e.~there is no $\varphi\in{\cal P}^{c,g}({\cal Z})$ which is equal to a
non--trivial convex linear combination of functions 
in ${\cal P}^{1,g}({\cal Z})$~.}
\begin{proof}
For decreasing functions, the property of being continuous from below is
obviously an extreme one.
\end{proof}
\theo{\Lemma}{23}{${\cal P}^{c,g}({\cal Z})$ is a Bauer simplex with respect
to the topology of full pointwise convergence, whose extreme points are
precisely the below--continuous functions of the form $\rho\circ g$ with
bounded characters $\rho$ on $T$~.}
\begin{proof}
Immediate by the above and Lemma 22~.
\end{proof}

In Lemma 23~, the below--continuous functions of the form $\rho\circ g$ with 
bounded characters $\rho$ on $T$ correspond to 
$\mu\in M^{1,e}_+({\cal Y})$ with
$Q_{Z_1},\dots,Q_{Z_m}$ \quad $\mu$--independent 
for any $Z_1,\dots,Z_m\in{\cal Z}$ with pairwise disjoint supports.
The question is whether this extreme boundary is not only closed with respect
to full pointwise convergence, but also with respect to weak convergence in 
$M^{1,e}_+({\cal Y})$~.  So let $\mu_1,\mu_2,\dots$ have the above property
and $\mu_n\to\mu$~.  Then, as {\it pars pro toto,}
$\mu_n(Q_{Z_1}\cap Q_{Z_2})=\mu_n(Q_{Z_1})\cdot\mu_n(Q_{Z_2})$~,
and in case that $\mu(Q_{Z_1}\cap Q_{Z_2})=\mu(Q_{Z_1})=\mu(Q_{Z_2})=0$
we get what we want.  By the same countability and approximation argument as
in $M^1_+(\R^d)$ this is enough to get the desired closedness property.\\
But for the uniqueness of mixtures of extreme points, there is still a problem
with the two different topologies involved, which give rise to different sets
of mixing Radon measures.

\begin{Acknowledgements}
I am most grateful to my PhD supervisor Prof Dr Paul Ressel (Eichst\"att) 
for having introduced me in a stimulating way into the powerful set of
methods which has led to [HR1999]~, [HR2000]~, [H2003]~, the present work and
a DFG application 
on {\rm invariant probability measures on combinatorial structures.}\\
I am also grateful to my present patron Prof Dr Joachim Gwinner 
(UniBw M\"unchen) for the work environment and the freedom 
which allowed me to write [H2003]~, the present paper and the mentioned DFG
application.
\end{Acknowledgements}

\begin{References}
\item {\sc Hirth U.}: Exchangeable random ordered trees by positive definite 
functions. {\it Journal of Theoretical Probability} {\bf 16 (2),} 339-344 
(April 2003).
\item {\sc Hirth U., Ressel P.}: Random partitions by semigroup methods. 
{\it Semigroup Forum} {\bf 59,} 126-140 (1999).
\item {\sc Hirth U., Ressel P.}: Exchangeable random orders and almost 
uniform distributions. {\it Journal of Theoretical Probability} 
{\bf 13,} 609-634 (2000).
\item {\sc Siegrist, K.}: Renewal processes on partially ordered sets.
{\it Random Structures and Algorithms} {\bf 22,} 15-32 (January 2003).
\end{References}
\begin{Address}
{\sc Dr.~rer.~nat. Ulrich Hirth},
Universit\"at der Bundeswehr M\"unchen, LRT 1, Werner--Heisenberg--Weg 39, 
85577 Neubiberg,
Germany; {\underline{\bf E--Mail}} {\tt Ulrich.Hirth@UniBw-Muenchen.de}~.
\end{Address}
\end{document}